\documentclass[11pt]{article}

\usepackage{amsmath}
\usepackage{amssymb, color, soul}
\usepackage{theorem}

\input xy
\xyoption {all}
\def\oM{\overline{\mathcal{M}}}
\def\cM{{\mathcal{M}}}

\def\C{\mathbb{C}}

\def\qed{{\hfill $\Diamond$}}
\def\Aut{{\rm Aut}}
\def\ch{{\rm ch}}
\def\End{{\rm End}}
\def\b1{{\mathbf 1}}
\def\rt{\rm{rt}}

\newtheorem{theorem}{Theorem}

\newtheorem{proposition}{Proposition}[section]

\newtheorem{lemma}[proposition]{Lemma}

{\theorembodyfont{\rmfamily}

\newtheorem{example}[proposition]{Example}
\newtheorem{remark}[proposition]{Remark}

}

\title{The Chern character of the Verlinde bundle over $\oM_{g,n}$}

\author{A. Marian
\thanks{Department of Mathematics, Northeastern University, a.marian@neu.edu}
\and
D. Oprea
\thanks{Department of Mathematics, University of California, San Diego, doprea@math.ucsd.edu}
\and
R. Pandharipande
\thanks {Department of Mathematics, ETH Z\"urich, rahul@math.ethz.ch}
\and
A. Pixton
\thanks{Department of Mathematics, Harvard University, apixton@math.harvard.edu}
\and
D. Zvonkine
\thanks{CNRS, Institut de Math\'ematiques de Jussieu, zvonkine@math.jussieu.fr}}
\date{}
\begin{document}
\baselineskip=16pt
\maketitle

\begin{abstract}
We prove an explicit formula for the total Chern character of the Verlinde bundle of conformal blocks over $\oM_{g,n}$ in terms of tautological classes.  The Chern characters of the Verlinde bundles define a semisimple CohFT (the ranks, given by the Verlinde formula, determine a semisimple fusion algebra). According to Teleman's classification of semisimple CohFTs, there exists an element of Givental's group 
 transforming the fusion algebra into the CohFT. We determine the element using the first Chern class of the Verlinde bundle on the interior~$\cM_{g,n}$ and the projective flatness of the Hitchin connection. 
\end{abstract}

\section{Introduction}

\subsection {The Verlinde bundle} Let $G$ be a complex, simple, simply connected Lie group. For any choice of genus $g$, and $n$ irreducible representations $\mu_1, \dots, \mu_n$ of the Lie algebra $\mathfrak g$ at level $\ell$, a vector bundle $${\mathbb E}_{g}(\mu_1, \dots, \mu_n)\to \oM_{g,n}$$ over the moduli space of stable pointed curves is constructed in~\cite {TUY}. Geometrically, over smooth pointed curves, the fibers of ${\mathbb E}_{g}(\mu_1, \dots, \mu_n)$ are the spaces of non-abelian theta functions -- spaces of global sections of the determinant line bundles over the moduli of parabolic $G$-bundles. The theory of {\it conformal blocks} is used in~\cite {TUY} to extend the bundle ${\mathbb E}_{g}(\mu_1, \dots, \mu_n)$ over the boundary of the moduli space $\oM_{g,n}$. With the geometric viewpoint in mind, we will refer to $\mathbb E_g(\mu_1, \ldots, \mu_n)$ as the Verlinde bundle. The same object is also known in the literature as the bundle of conformal blocks or the bundle of vacua.

Two basic invariants of the Verlinde bundles have been studied so far. First, $$\text{ rank } \mathbb E_{g} (\mu_1, \ldots, \mu_n)=d_g(\mu_1, \ldots, \mu_n)$$ is given by the Verlinde formula; see for instance \cite {beauville}. Second, an explicit closed expression for the first Chern class was obtained in genus $0$ in \cite{F},~\cite{Mu}, in genus $1$ and for one marking in~\cite{F}, and in arbitrary higher genus in~\cite {MOP}. We calculate here the total Chern character. As a consequence of our formulas, the higher Chern classes are seen to lie in the tautological ring
$$RH^\star(\oM_{g,n}) \subset H^\star(\oM_{g,n})\, , $$ see~\cite{FP}.
We expect our formulas to hold in the Chow ring as well; however our proof does not yield this, since it relies on the classification of semisimple CohFTs, which is only established in cohomology.

The expression for the Chern character is derived as follows. It is well known that the restriction of the Verlinde bundle to each boundary divisor of $\oM_{g,n}$ decomposes into a direct sum of tensor products of analogous Verlinde bundles \cite {TUY}. These factorization rules imply that the total Chern characters 
$$
\Omega_{g, n}(\mu_1, \dots, \mu_n) = \ch \,{ \mathbb E}_{g}(\mu_1, \dots, \mu_n)
$$
define a cohomological field theory (CohFT). Furthermore, the
CohFT obtained
is  semisimple. To explicitly solve the theory, we make use of the following three ingredients.
\begin {itemize}
\item [(i)] The first Chern class of the Verlinde bundle over the moduli of smooth pointed curves $\mathcal M_{g, n}$ is, in the form we need, a consequence of \cite {T}. 
\item [(ii)] The existence of a projectively flat connection in the Verlinde bundle over $\mathcal M_{g, n}$ \cite {TUY} allows one to compute the full Chern character over $\mathcal M_{g, n}$.
\item [(iii)] Finally, the Givental-Teleman classification of semisimple CohFTs ~\cite{Teleman} {\it uniquely} determines the extension of the Chern character to $\oM_{g,n}$ (the uniqueness is the core of Teleman's paper).
\end {itemize}

\subsection {The main result} To state the result, we recall the following representation theoretic quantities. First, for a simple Lie algebra $\mathfrak g$ and for each level $\ell$, the conformal anomaly is given by $$c=c(\mathfrak g, \ell)=\frac{\ell \dim \mathfrak g}{\check{h}+\ell}$$ where $\check{h}$ is the dual Coxeter number. Next, for each  representation with highest weight $\mu$ of level $\ell$, we set
$$
\mathsf w(\mu)=\frac{(\mu, \mu+2\rho)}{2(\check{h}+\ell)}.$$ Here $\rho$ is half of the sum of the positive roots, and the Cartan-Killing form $(, )$ is normalized so that the longest root $\theta$ has $$(\theta, \theta)=2.$$ 

\begin{example} For ${\mathfrak g} = {\mathfrak s}{\mathfrak l} (r, {\mathbb C}),$ the highest weight of a representation of level $\ell$ is given by an $r$-tuple of integers 
$$\mu = (\mu^{1}, \ldots, \mu^{r}), \, \, \ell \geq \mu^{1} \geq \cdots \geq \mu^{r} \geq 0,$$ defined up to shifting the vector components by the same integer. 
Furthermore, we have
$$ c (\mathfrak g, \ell) = \frac{\ell (r^2-1)}{\ell+r},$$
$${\mathsf w} (\mu)=\frac{1}{2(\ell + r)}\left(\sum_{i=1}^{r} (\mu^i)^2-\frac{1}{r}\left(\sum_{i=1}^r \mu^i\right)^2+ \sum_{i=1}^{r} (r-2i+1)\mu^i\right).$$
\end{example}

\medskip

The Verlinde Chern character will be expressed as a sum over stable graphs $\Gamma$ of genus~$g$ with $n$ legs. The contribution of each graph is a product of several factors coming from legs, edges and vertices, which we now describe. 

\begin{itemize}
\item[(i)] {\bf Leg factors.} Each leg $l$ determines a marking, and is assigned
\begin {itemize} 
\item a cotangent class $\psi_l$;
\item  a representation $\mu_l$ of the Lie algebra $\mathfrak g$. 
\end {itemize} The corresponding factor is
$$
\mathsf {Cont }(l) = \exp\left(-{\mathsf w(\mu_l) \cdot \psi_l}\right).
$$

\item[(ii)] {\bf Vertex factors.} Every vertex~$v$ in a stable graph comes with 
\begin {itemize}
\item a genus assignment $g_v$;
\item $n_v$ adjacent half-edges including the legs labeled with representations
$\mu_1, \dots, \mu_{n_v}$.
\end {itemize}
The vertex factor is defined as the rank of the Verlinde bundle
$$
\mathsf{Cont} (v) = d_{g_v}(\mu_1, \dots, \mu_{n_v}).
$$
\item[(iii)] {\bf Edge factors.} To each edge~$e$ we assign:
\begin {itemize}
\item a pair of conjugate representations of level $\ell$: $\mu$ at one half-edge and the dual representation $\mu^*$ at the other half-edge; 
\item each edge corresponds to a node of the domain curve, thus determining two cotangent classes $\psi_e'$ and $\psi_{e}''$. 
\end {itemize} The edge factor is then defined as 
$$
\mathsf {Cont } (e) = \frac{1 - \exp (-{\mathsf w \cdot (\psi_e' + \psi_e'')})}{\psi_e' + \psi_e''},
$$
where $\mathsf w = \mathsf w(\mu) = \mathsf w(\mu^*)$. The edge factor depends on the representation assignment $\mu$, but this will not be indicated explicitly in the notation. 

\end{itemize}
Next, for each stable graph $\Gamma$ of type $(g, n)$, we consider all possible assignments $$\mu:E\to P_{\ell}$$ of conjugate representations of level $\ell$ to the edges $E$ of the graph $\Gamma$. For each such graph $\Gamma$ and representation assignment $\mu$, we form the cohomology class $$(\iota_{\Gamma})_{\star}\left(\prod_{l} \mathsf {Cont }(l)\prod_{v} \mathsf {Cont }(v) \prod_{e} \mathsf {Cont} (e)\right)\in H^{\star}(\oM_{g,n}),$$ where $\iota_{\Gamma}$ denotes the inclusion of the boundary stratum of curves of type $\Gamma$. 

Let $\lambda_1 \in H^2(\oM_{g,n})$ denote the first Chern class of the Hodge bundle
$$\mathbb{H}_g \rightarrow \oM_{g,n}$$
with fiber $H^0(C,\omega_C)$ over the moduli point $[C,p_1,\ldots,p_n] \in \oM_{g,n}$.

\begin{theorem} \label{Thm:main}
The Chern character of the Verlinde bundle is given by
\begin{multline*}
\ch\, {\mathbb E}_{g}(\mu_1, \dots, \mu_n) =\\
\exp \left(-\frac{c (\mathfrak g, \ell)}{2}\cdot \lambda_1\right)\sum_{\Gamma, \mu} \frac1{|\Aut (\Gamma)|} (\iota_{\Gamma})_{\star}\left(\prod_{l} \mathsf {Cont }(l)\prod_{v} \mathsf {Cont }(v) \prod_{e} \mathsf {Cont} (e)\right).
\end{multline*}
\end{theorem}
\vskip.1in

\begin{example} \label{Ex:sl2level1}
We illustrate the formula of the Theorem in arbitrary genus $g$, for the Lie algebra $\mathfrak g=\mathfrak{sl}_2,$ at level $\ell=1$. In this case, there are two representations to consider, denoted $\emptyset$ and $\square$. Then
$$
c(\mathfrak{sl}_2,1) = 1, \quad \mathsf w(\emptyset)=0,\quad \mathsf w(\square)=\frac{1}{4}.
$$ We compute
$$
\ch\, {\mathbb E}_g(\mu_1, \dots, \mu_n), 
$$
where each $\mu_i$ equals either $\emptyset$ or $\square$.
To this end, we remark that:
\begin {itemize} 
\item [--] No edge can be decorated by the representation $\emptyset$ since formula (iii) above yields a zero contribution. Thus, all edges are decorated by the representation $\square$.
\item [--] The degree of every vertex should become even after removing all the legs (markings) decorated by $\emptyset$. Indeed, the rank of the Verlinde bundle vanishes if the number of $\square$-representations is odd. A graph whose vertices satisfy this parity condition will be called $\square$-even.
\item [--] 
Furthermore, as soon as the number of $\square$-representations is even, the rank of the Verlinde bundle on a genus~$g$ moduli space equals $2^g$, independently of the number of marked points. This can be easily deduced from the degeneration formula, or via the expressions in~\cite{beauville}. The product of $2^{g_v}$ over the vertices of a stable graph~$\Gamma$ yields $2^{g-h^1(\Gamma)}$, where $h^1$ is the first Betti number.
\end{itemize} 
As a consequence, over $\oM_{g, n}$ we obtain
$$
\text{ch }\mathbb E_g=
\exp\left(-\frac{\lambda_1}{2}\right)\cdot 
\sum_{\Gamma}
\frac{2^{g-h^1(\Gamma)}}{|\text{Aut}(\Gamma)|}\cdot 
(\iota_{\Gamma})_{\star} 
\left(\prod_e \frac{1-\exp\left(-\frac{1}{4}(\psi'_e+\psi''_e)\right)}{\psi_e'+\psi''_e}\cdot \prod_l e^{-\psi_l/4} \right),
$$ 
where the sum is taken over the $\square$-even stable graphs, the first product runs over all edges, and the second product runs over the $\square$-marked legs. 

This example will be used in Section \ref{divsubr} where we examine whether the Chern character belongs to the cohomology subring generated by divisors, with a negative answer.

\end{example}

Another example -- that of the Verlinde bundle for $\mathfrak{sl}_r$ at level~1 over moduli space of curves of compact type -- will be worked out in Proposition \ref{exponential}.

\subsection {Acknowledgements} The calculation of the Chern character 
was completed in October 2013 during the workshop {\it Cohomology of the moduli space of curves} organized by the {\em Forschungsinstitut 
f\"ur Mathematik} at ETH Z\"urich. Additional funding for the
workshop was provided by the Clay Foundation.
We thank P. Belkale and N. Fakhruddin for discussions related to
the first draft of the paper.

A.M. was supported by the Sloan Foundation and the NSF through grant DMS 1303389. D. O. was supported by the Sloan Foundation and the NSF through grants DMS 1001486 and DMS 1150675. 
R.P. was supported by the Swiss National Science Foundation and
the European Research Council through
grants SNF-200021-143274 and ERC-2012-AdG-320368-MCSK.
A.P. was supported by the Clay Foundation.
D.~Z. was supported by the grant ANR-09-JCJC-0104-01.

\section{Semisimple CohFTs}

\subsection {Definitions} We succinctly recall terminology related to cohomological field theories. Fix a finite dimensional complex vector space $V$, endowed with a non-degenerate pairing $\eta$ and a distinguished element $1\in V$. A 
{\em CohFT} is the data of cohomology classes $$\Omega=(\Omega_{g, n}),\,\,\,\, \Omega_{g, n}\in H^{\star}(\oM_{g, n})\otimes (V^{\star})^{\otimes n}$$ for $2g-2+n>0$, subject to the requirements:
\begin {itemize}
\item [(i)] each $\Omega_{g, n}$ is invariant under the action of the symmetric group $S_n$;
\item [(ii)] $\Omega$ is compatible with the gluing maps. Explicitly, for the gluing map $$gl:\oM_{g-1, n+2}\to \oM_{g, n},$$ the pullback ${gl}^{\star} \Omega_{g, n}$ equals the contraction of $\Omega_{g-1, n+2}$ with $\eta^{-1}$ at the two extra-markings. The same requirement is enforced for the second gluing map $$gl:\oM_{g_1, n_1+1}\times \oM_{g_2, n_2+1}\to \oM_{g, n}.$$ The tensor $\eta^{-1}\in V\otimes V$ above is given in a fixed basis $\{e_{\mu}\}$ of $V$ by the expression $$\eta^{-1}=\sum_{\mu, \nu} \eta^{\mu\nu} e_{\mu}\otimes e_{\nu}$$ where $(\eta^{\mu\nu})$ is the inverse of the matrix $\eta_{\mu\nu}=\eta(e_{\mu}, e_{\nu}).$ Consequently, the first gluing rule rewrites as $$gl^{\star}\Omega_{g, n}(v_1\otimes \ldots\otimes v_n)=\sum_{\mu, \nu} \eta^{\mu\nu}\cdot \Omega_{g-1, n+2}(v_1\otimes \ldots\otimes v_n\otimes e_{\mu}\otimes e_{\nu}).$$ The second gluing formula takes a similar shape.
\item [(iii)] for all $v_1, \ldots, v_n\in V$, we have $$\Omega_{g, n+1} (v_1\otimes \ldots \otimes v_n\otimes 1)=p^{\star} \Omega_{g, n}(v_1\otimes \ldots \otimes v_n),$$ where $p: \oM_{g, n+1} \to \oM_{g,n}$ is the forgetful map, and $$\Omega_{0, 3}(v_1\otimes v_2\otimes 1)=\eta(v_1, v_2).$$
\end {itemize}

When $\Omega_{g, n}$ are cohomology classes of degree $0$, the definition above recovers the axioms of topological quantum field theory (TQFT). 

Each CohFT defines on $V$ the structure of an associative algebra with unit via the quantum product $$\eta(v_1\bullet v_2, v_3)=\Omega_{0, 3}(v_1\otimes v_2\otimes v_3).$$ Associativity follows from applying~(ii) to two gluing maps $\oM_{0,3} \times \oM_{0,3} \to \oM_{0,4}$, while (iii) implies the existence of unit. We are concerned with theories for which the algebra $V$ is semisimple. 

Such theories are classified in \cite{Teleman}. Specifically, $\Omega$ is obtained from the algebra $V$ via the action of an $R$-matrix $$R \in {\bf 1} + z \cdot \End V [[z]],$$ satisfying the {\it symplectic condition} $$R(z)R^{\star}(-z)=\bf 1.$$ Here $R^{\star}$ denotes the adjoint with respect to $\eta$ and $\bf 1$ is the identity matrix. This condition is called {\it symplectic} because it is equivalent to requiring that the multiplication by $R(z)$ preserve a natural symplectic form in $V[[z]][z^{-1}]$ (see~\cite{Givental}). It can also be written as 
$$
R(z) \eta^{-1} R(z)^{\top} = \eta^{-1},
$$
where $\top$ indicates transpose, and thus $R(z) \eta^{-1} R(z)^{\top}$ is the action of $R(z) \otimes R(z)$ on the bivector $\eta^{-1} \in V \otimes V$.

The explicit reconstruction of the semisimple CohFT from the $R$-matrix action will be explained below, following \cite {PaPiZv}.

\subsection {Actions on CohFTs}

We begin by describing two basic actions on cohomological field theories. 

Assume that $\Omega=(\Omega_{g, n})$ is a CohFT with underlying vector space $(V, 1, \eta)$.
Fix a symplectic matrix $R\in {\bf 1}+z\cdot \text{End }(V)[[z]]$ as above. A new CohFT on the space $(V, 1, \eta)$ is obtained via the cohomology elements $$R \Omega=(R\Omega)_{g, n},$$  defined as sums over stable graphs $\Gamma$ of genus $g$ with $n$ legs, with contributions coming from vertices, edges and legs. Specifically, \begin{equation}\label{romega}(R\Omega)_{g, n}=\sum_{\Gamma} \frac{1}{|\text{Aut }(\Gamma)|} (\iota_{\Gamma})_{\star} \left(\prod_{l}\mathsf {Cont}(l)\prod_{v}\mathsf {Cont}(v)\prod_{e}\mathsf {Cont }(e) \right)\end{equation} where:
\begin {itemize}
\item [(i)] the vertex contribution is $$\mathsf {Cont}(v)=\Omega_{g(v), n(v)},$$ with $g(v)$ and $n(v)$ denoting the genus and number of half-edges and legs of the vertex;
\item [(ii)] the leg contribution is the $\text{End}(V)$-valued cohomology class $$\mathsf{Cont}(l)=R(\psi_l)$$ where $\psi_l$ is the cotangent class at the marking corresponding to the leg; 
\item [(iii)] the edge contribution is $${\mathsf {Cont}}(e)=\frac{\eta^{-1}-R(\psi'_e)\eta^{-1} R(\psi''_e)^{\top}}{\psi'_e+\psi''_e}.$$ Here $\psi'_e$ and $\psi''_e$ are the cotangent classes at the node which represents the edge $e$. The symplectic condition guarantees that the edge contribution is well-defined. 

For the benefit of the reader, we clarify the meaning of the expression 
$$
\mathsf{Cont}(e)\in V^{\otimes 2}\otimes H^{\star}(\oM_{g', n'})\otimes H^{\star}(\oM_{g'', n''}),
$$ 
where $(g', n')$ and $(g'', n'')$ are the labels of the vertices adjacent to~$e$. To this end, we describe the edge contribution in coordinates, after fixing a basis $\{e_\mu\}$ of $V$. The components of the $R$-matrix in this basis are $R_\mu^\nu(z)$, in other words,
$$
R(z)(e_\mu)=\sum_{\nu} R_\mu^\nu (z) \cdot e_\nu.
$$ 
Then the components of $\mathsf{Cont}(e)$ are given by
$$
\mathsf{Cont}(e)^{\mu \nu} = 
\frac{\eta^{\mu \nu}- \sum_{\rho, \sigma} R_\rho^\mu(\psi_e')\cdot \eta^{\rho \sigma}\cdot R_\sigma^\nu(\psi_e'')}{\psi_e'+\psi_e''}\in H^{\star}(\oM_{g', n'})\otimes H^{\star}(\oM_{g'', n''}).
$$ 
The fraction 
$$
\frac{\eta^{\mu \nu}-\sum_{\rho, \sigma}R_\rho^\mu(z)\cdot \eta^{\rho \sigma}\cdot R_\sigma^\nu(w)}{z+w}
$$ 
is a power series in $z$ and $w$, since the numerator vanishes as $z=-w$. This is a consequence of the symplectic condition which in coordinates is seen to take the form 
$$
\sum_{\rho, \sigma} R_\rho^\mu(z)\cdot \eta^{\rho \sigma}\cdot R_\sigma^\nu(w) = \eta^{\mu \nu}.
$$ 
The substitution $z=\psi_e'$ and $w=\psi_e''$ is therefore unambiguously defined.
\end {itemize}

\begin{remark}
To simplify our formulas, we have changed Givental's and Teleman's conventions by replacing $R$ with $R^{-1}$. In particular, equation \eqref{romega} above determines a right group action on CohFTs, rather than a left group action as in Givental's and Teleman's papers. This will play no role here.
\end{remark}

A second action on CohFTs is given by translations. As before, let $(\Omega, V, 1, \eta)$ be a CohFT, and consider a power series $T\in V[[z]]$ with no terms of degrees $0$ and $1$: $$T(z)=T_2z^2+T_3z^3+\ldots,\,\,\,T_k\in V.$$ A new CohFT based on $(V, 1, \eta)$, denoted $T\Omega$, is defined by setting \begin{equation}\label{translation} (T\Omega)_{g, n} (v_1\otimes \ldots \otimes v_n)=\sum_{m=0}^{\infty} \frac{1}{m!} (p_m)_{\star} \Omega_{g, n+m}(v_1\otimes \ldots\otimes v_n \otimes T(\psi_{n+1})\otimes \ldots \otimes T(\psi_{n+m}))\end{equation} where $$p_{m}:\oM_{g, n+m}\to \oM_{g, n}$$ is the forgetful morphism. Expression~\eqref{translation} should be understood as formal expansion, by distributing the powers of $\psi$-classes as follows: $$\Omega_{g, n+ m}(\cdots \otimes T(\psi_{\bullet})\otimes\cdots )=\sum_{k=2}^{\infty} \psi_{\bullet}^k \cdot \Omega_{g, n+m} (\cdots \otimes T_k \otimes \cdots).$$
 
 \subsection {Reconstruction} With the above terminology understood, we can state the Givental-Teleman classification theorem. 

\begin{theorem}[\cite {Teleman}]
Let $\Omega$ be a semisimple CohFT on $(V, 1, \eta)$. Then there exists a symplectic matrix 
$$
R\in {\mathbf 1}+z\cdot \End (V)[[z]]
$$
such that
$$
\Omega=RT\omega,
$$
where $\omega$ is the degree $0$ topological part of $\Omega$ and the power series $T$ is obtained by evaluating $R(z)$ at $1 \in V$, removing the free term, and multiplying by~$z$:
$$
T(z)=z(1-R(z)\cdot 1)\in V[[z]].
$$ 
\end{theorem}
 
In view of this result, we prove:
 \begin {lemma} \label{recon} The $R$-matrix of a semisimple CohFT $\Omega$ is uniquely determined by the restrictions of the elements $\Omega_{g, n}$ to $\cM_{g, n}.$
\end {lemma} 

\noindent {\bf Proof.} We establish the lemma using the reconstruction statement $$\Omega=RT\omega,$$ and working in an idempotent basis for $V,$ afforded by semisimplicity. Explicitly, let $v_\alpha$ be so that $$v_{\alpha}\bullet v_{\beta}=\delta_{\alpha, \beta} v_{\alpha}.$$ Denote by
$c_\alpha = \eta(v_\alpha, v_\alpha)$ the scalar square of~$v_\alpha$. Note that $c_\alpha \ne 0$.
The TQFT axioms imply 
$$
\omega_{g,n}(v_\alpha, v_\beta, \dots, v_\beta)=c_\beta^{1-g} \, \delta_{\alpha, \beta}.
$$ 
This is seen by applying the degeneration rule (ii) in the definition of a CohFT to any of the gluing maps $(\oM_{0,3})^{2g-2+n} \to \oM_{g,n}$. We expand $$R (z) (v_{\alpha}) =\sum_{\beta, k} (R_k)_{\alpha}^{\beta}\cdot v_{\beta}\cdot z^k.$$ Clearly, it suffices to explain the uniqueness of $(R_k)_{\alpha}^{\beta}$. Consider the cohomology class
$$
\Omega_{g,n}(v_\alpha, v_\beta, \dots, v_\beta) \in H^*(\cM_{g,n})
$$ over the open part of the moduli space. 
If $g > 3k$ the coefficient of $\psi_1^k$ in this cohomology class is well-defined, because the tautological ring $R^{\star}(\cM_{g,n})$ is generated by the classes $$\kappa_1, \ldots, \kappa_{[g/3]}, \psi_1, \ldots, \psi_n$$ with no relations up to degree $g/3$. The freeness up to degree $g/3$ is a consequence of the
stability results of \cite{bo,l}, see also \cite {BF}. The exact expression for the Givental group action $$\Omega=RT\omega$$ given in \eqref{romega} and \eqref{translation} is used to find the coefficient of $\psi_1^{k}$ in the restriction of $\Omega$ to $\cM_{g, n}$. We claim this coefficient equals
$c_\beta^{1-g} \cdot (R_k)_\alpha^\beta$. This follows from the following observations:
\begin {itemize}
\item [--] in expression \eqref{romega}, the only graph $\Gamma$ contributing to the restriction of $$\Omega=RT\omega$$ to $\cM_{g, n}$ is the single vertex graph;
\item [--] in the translation action \eqref{translation}, the terms $m\geq 1$ contribute pushforwards of monomials in the $\psi$-classes $$(p_{m})_{\star}(\psi_{n+1}^{\ell_1}\cdots \psi_{n+m}^{\ell_m}).$$ These can be expressed as polynomials in the $\kappa$-classes and thus do not involve $\psi_1^k$, see \cite {AC} for explicit formulas;
\item [--] finally, for the terms corresponding to the single vertex graph and $m=0$, we extract the coefficient of $\psi_1^k$ with the aid of the identity $$\omega_{g,n}(v_\alpha, v_\beta, \dots, v_\beta)=c_\beta^{1-g} \, \delta_{\alpha, \beta}.
$$
\end {itemize}
As a consequence, every $(R_k)_\alpha^\beta$ is uniquely determined by the restriction of the classes $\Omega_{g,n}$ to $\cM_{g,n}$ for $g$ large enough, as claimed. 
\qed

 \section{Proof of the Theorem}
 
 \subsection {The CohFT obtained from conformal blocks} The total Chern character of the bundle of conformal blocks defines a CohFT.  Explicitly:
\begin {itemize} 
\item[(a)] the vector space $V$ has as basis the irreducible representations of $\mathfrak g$ at level $\ell$. The distinguished element ${1}$ corresponds to the trivial representation. The pairing $\eta$ is given by $$\eta(\mu, \nu)=\delta_{\mu, \nu^{\star}}$$ where $\nu^{\star}$ denotes the dual representation. 
\item [(b)] The cohomology elements defining the theory are $$\Omega_{g, n}(\mu_1, \ldots, \mu_{n})=\ch_t \, (\mathbb E_{g}(\mu_1, \ldots, \mu_n))\in H^{\star}(\oM_{g, n}).$$ Here, for a vector bundle $\mathbb E$ with Chern roots $r_1, \ldots, r_k$ we write $$\ch_t(\mathbb E)=\sum_{j=1}^{k} \exp(t\cdot r_j).$$ 
The parameter $t$ can be treated either as a formal variable, in which case we work over the ring $\C[[t]]$ instead of~$\C$, or as a complex number, so that all the statements are true for any $t \in \C$.
\end {itemize}
Axiom (i) in the definition of CohFT is obvious. Axiom (ii) follows from the fusion rules of \cite {TUY}. For instance, for the irreducible boundary divisor we have $$gl^{\star} \mathbb E_g(\mu_1, \ldots, \mu_n)=\bigoplus_{\nu\in P_{\ell}} \mathbb E_{g-1}(\mu_1, \ldots, \mu_n, \nu, \nu^{\star})$$ 
and taking Chern characters we find $$\Omega_{g, n}(\mu_1\otimes \ldots\otimes \mu_n)=\sum_{\nu\in P_{\ell}} \Omega_{g-1, n+2}(\mu_1\otimes \ldots\otimes \mu_n\otimes \nu\otimes \nu^{\star}),$$ as required by (ii). 
The existence of a unit as required by axiom (iii)
 is {\em propagation of vacua}, and is proved in the form needed here in \cite {F}, Proposition 2.4(i). Indeed, under the marking-forgetting map $$p:\oM_{g, n+1}\to \oM_{g, n}\, ,$$ we have $$p^{\star} \mathbb E_g(\mu_1, \ldots, \mu_n)=\mathbb E_g(\mu_1, \ldots, \mu_n, {1}).$$ Finally, the requirement that $$\Omega_{0, 3}(\mu\otimes \nu \otimes {1})=\delta_{\mu, \nu^{\star}}$$ is Corollary $4.4$ of \cite {beauville}.  
  
The CohFT $\Omega$ thus constructed is semisimple. Indeed, taking the degree~0 part of~$\Omega$ (or, equivalently, setting $t=0$, since the powers of~$t$ track the cohomological degree) we get the TQFT 
$$
\omega_{g,n}=\Omega_{g, n}|_{t=0}
$$ 
given by 
$$
\omega_{g, n}\in (V^{\star})^n,\,\,\, \omega_{g, n}(\mu_1, \ldots, \mu_n)=\text{rk } \mathbb E_g(\mu_1, \ldots, \mu_n)=d_g(\mu_1, \ldots, \mu_n).
$$ 
The associated Frobenius algebra is the Verlinde fusion algebra, and is known to be semisimple. An account can be found in Proposition $6.1$ of \cite {beauville}.

\subsection {The R-matrix of the CohFT of conformal blocks}
Teleman's classification~\cite{Teleman} ensures that $\Omega$ is obtained from the Verlinde fusion algebra $\omega$ by Givental's group action of an $R$-matrix $R \in \End\, V [[z]]$. The $R$-matrix of the theory will be found below, and shown to be diagonal in the natural basis of $V$ consisting of irreducible representations at level $\ell.$

It will be more convenient to consider a slightly modified CohFT $\Omega'$, given by $$\Omega'_{g, n}=\Omega_{g, n}\exp\left(t\cdot \frac{c(\mathfrak g, \ell)}{2}\cdot \lambda_1\right).$$ The fact that  $\Omega'$ still satisfies the requisite axioms follows from the fact that the Hodge bundle splits compatibly over the boundary divisors in $\oM_{g, n}$. Indeed, only axiom (ii) requires checking. For instance, compatibility of $\Omega'$ with the gluing map $$gl:\oM_{g-1, n+2}\to \oM_{g, n}$$ is tantamount to the matching of prefactors $$gl^{\star}\exp\left(t\cdot \frac{c(\mathfrak g, \ell)}{2}\cdot \lambda_1\right)=\exp\left(t\cdot \frac{c(\mathfrak g, \ell)}{2}\cdot \lambda_1\right),$$ which is in turn justified by the identity $$gl^{\star} \lambda_1=\lambda_1;$$ see \cite {AC} for instance. The case of the second gluing map is entirely similar.

\vskip.1in

We proceed to find the restrictions of the theory $\Omega'$ to $\cM_{g,n}$. The Verlinde bundle is projectively flat over $\cM_{g,n}$ by results of \cite {TUY}, \cite {T}. Therefore over $\cM_{g,n}$ its Chern character is given by 
\begin{equation} \label{Eq:flat}
\ch_t(\mathbb E_g(\mu_1, \ldots, \mu_n))=d_{g}(\mu_1, \dots, \mu_n) \cdot \exp \left(t\cdot \frac{c_1(\mathbb E_g(\mu_1, \ldots, \mu_n))}{d_g(\mu_1, \ldots, \mu_n)}\right).
\end{equation} (The Chern character is set to zero when the denominator $d_g(\mu_1, \ldots, \mu_n)$ of the fraction above vanishes.) Formula \eqref{Eq:flat} is standard and is explained for instance in Chapter II.3 of \cite {K}. Now, over the open part of the moduli space $\cM_{g,n}$, the slope of the Verlinde bundle was written in \cite {MOP} as a consequence of \cite {T}: $$\frac{c_1(\mathbb E_g(\mu_1, \ldots, \mu_n))}{d_g(\mu_1, \ldots, \mu_n)}=-\frac{c(\mathfrak g, \ell)}{2}\cdot \lambda_1-\sum_{i=1}^{n} \mathsf w(\mu_i)\psi_i.$$ The signs differ from \cite {MOP}, where the dual bundle of covacua was used. Therefore, \begin{equation}\label{omega'}\Omega'_{g, n}(\mu_1\otimes \ldots\otimes \mu_n)=d_g(\mu_1, \ldots, \mu_n)\cdot \exp \left(-t\cdot \sum_{i=1}^{n} \mathsf w(\mu_i)\psi_i\right)\end{equation} over $\cM_{g, n}$. 
\vskip.1in

We must have $\Omega'=RT\omega$ for a symplectic matrix $R$ and the translation $T$ introduced previously. Let $W \in \text{End}\, V[[z]] $ be the diagonal matrix in the basis of level $\ell$ representations, whose diagonal elements are given by $$W(z)_{\mu}^\mu=\exp(-tz\cdot \mathsf w(\mu)).$$ Since $\mathsf w(0) = 0$, we see that $W(z)\cdot 1 = 1$.  
Hence, the associated $T(z)$ vanishes,
$$T(z)=z(1-W(z)\cdot 1)=0\, .$$
By the above discussion, the identification of the two CohFT's $$\Omega'=W\omega$$ holds over $\cM_{g, n}$ for all $g$ and $n$. This is precisely the content of equation \eqref{omega'} for the left-hand side, and equation \eqref{romega} for the right-hand side. By the unique reconstruction of the $R$-matrix from restrictions, proved in Lemma \ref{recon}, we conclude that $R=W$.
 Hence the equality $$\Omega'=W\omega$$ also holds over the compact moduli space $\oM_{g,n}$. Formula \eqref{romega} applied to the matrix $W$ and to the theory $\omega$ yields the expression of the Theorem.

\begin{remark}
It is, of course, possible to find the $R$-matrix that takes the Verlinde fusion algebra directly to the CohFT  $\Omega$ rather than~$\Omega'$. This $R$-matrix is given by
$$
{R(z)}_\mu^\mu = 
\exp \left(t z \cdot 
\left(-{\mathsf w}(\mu) + \frac{c({\mathfrak g}, \ell)}{24}\right) \right).
$$
\end{remark}

\begin{remark}
The Verlinde slope formula over ${\mathcal M}_{g,n}$, written in ~\cite{MOP}, is used as input for our derivation. The full slope formula in \cite{MOP}, on the compactification $\oM_{g,n}$, is then recovered by the result proven here. This can be seen by explicitly accounting for all the one-edge stable graphs $\Gamma$ that contribute to the equation of Theorem \ref{Thm:main}.
The matching of the formula here with \cite{MOP} provides a nontrivial
check.
\end{remark}

\begin{remark}
Since the projective flatness of the Hitchin connection is used as an input, 
our formula for the Chern character of the Verlinde bundle does not immediately yield nontrivial relations in
the tautological ring $RH^*({\mathcal M}_{g,n})$ by imposing projective flatness.
 \end{remark}

\section{The divisor subring}\label{divsubr}

If a vector bundle $V$ on a nonsingular algebraic variety $X$ admits a flat connection with logarithmic singularities along 
divisors 
$$\bigcup_i D_i\subset X$$ with simple normal crossings, then by a theorem of
Esnault and Verdier~\cite{EsnVie}, the Chern character of $V$ lies in the 
subring of $H^\star(X)$ generated by the divisors $D_i$.  
Fakhruddin \cite{F} applies the result 
to express the Chern character of the Verlinde bundle in genus~0 in terms of the residues of the Hitchin connection. Furthermore, Fakhruddin claims
the result of Esnault-Verdier implies that the Chern character of the
Verlinde bundle lies in the subring of $H^\star(\overline{\cM}_{g,n},\mathbb{Q})$
generated by divisors  for {\em all} $g$ and $n$
(see Question~7.3 of \cite{F} and the preceding paragraph). 

We will explain that as a consequence of Theorem \ref{Thm:main} the Chern character of the Verlinde bundle does {\em not} lie in the subring of $H^\star(\overline{\cM}_{g,n})$
generated by divisors. There is no surprise here: since certain boundary divisors of $\oM_{g,n}$ (for instance, the divisor $\delta_0$ of stable curves with a non-separating node) have self-intersections, the conditions of the Esnault-Verdier theorem are {\em not} satisfied (and the implication in the paragraph before Question~7.3 in \cite{F} is incorrect).

We prove first a positive result. Let $\cM_{g,n}^{\rt} \subset \overline{\cM}_{g,n}$ be the moduli space of curves with rational tails, that is curves with a genus~$g$ irreducible component. We show:

\begin{proposition} \label{Prop:CompType}
The restriction of $\ch (\mathbb{E}_{g,n}(\mu_1, \dots, \mu_n))$ to $\cM_{g,n}^{\rt}$ lies in the subring of $H^\star(\cM_{g,n}^{\rt})$ generated by divisors.
\end{proposition}

\paragraph{Proof.} The Chern character $\ch (\mathbb{E}_{g,n}(\mu_1, \dots, \mu_n))$ over $\cM_{g,n}^{\rt}$ is computed by restricting the summation in Theorem \ref{Thm:main} from all stable graphs to stable trees possessing one vertex of genus $g$. Let $\Gamma$ be such a tree, and consider the class \begin{equation} \label{Eq:term}
(\iota_\Gamma)_{\star}\left(\prod_{i=1}^n \psi_i^{d_i} \prod_e (-\psi'_e-\psi''_e)^{k_e}\right).
\end{equation}
This can be represented as the product
\begin{equation} \label{Eq:prod}
\prod_{i=1}^n \psi_i^{d_i} \prod_e \delta_e^{k_e+1},
\end{equation}
where $\delta_e$ is the boundary divisor in $\cM_{g, n}^{\rt}$ corresponding to the 1-edge stable graph obtained from $\Gamma$ by contracting all edges except for $e$. The key identity $$\delta_e^{k_e+1}=\iota_{\star}(-\psi'_e-\psi''_e)^{k_e}$$ holds over $\cM_{g,n}^{\rt}$. Indeed, this formula may fail when $e$ carries an unmarked vertex of genus $\leq g/2$; in this case, there are corrections coming from strata contained in the self-intersection of $\delta_e$. Since we restrict to curves with rational tails, such unmarked vertices are not allowed: the boundary of $\cM_{g,n}^{\rt}$ has simple normal crossings. 
Now it suffices to note that the Chern character is explicitly a polynomial in $\lambda_1$ multiplied by a linear combination of terms of the form~\eqref{Eq:term}. Equation \eqref{Eq:prod} shows that each term is a product of divisor classes. \qed

\begin{proposition} \label{Prop:NotDiv}
The second Chern character of the Verlinde bundle of $\mathfrak{sl}_2$ at level~1 in general does not lie in the subring of $H^*(\oM_{g,n})$ generated by divisors.
\end{proposition}

\paragraph{Proof.} The formula for the Chern character of $\mathfrak{sl}_2$ at level~1 was worked out in Example~\ref{Ex:sl2level1}. We will only use the case where all representations $\mu_i$ are equal to~$\Box$; therefore we will assume that $n$ is even. 

In our expression for the Chern character, 
we will ignore the exponential factor containing $\lambda_1$ since both this factor and its inverse lie in 
the subring generated by divisors.

Consider the stable graphs with exactly two vertices joined by exactly two edges. The vertices can be of any genus and support any number of legs. Such graphs will be called 2-{\em loops}. If a 2-loop has an even (respectively, odd) number of legs attached to both vertices we will call it {\em even} (respectively, {\em odd\/}). According to Example~\ref{Ex:sl2level1}, the coefficient of odd 2-loops  in the expression of $\ch_2(\mathbb{E}_{g,n}(\Box, \dots, \Box))$ vanishes, while the coefficient of even 2-loops equals~$1/16$. We will show that this is incompatible with $\ch_2(\mathbb{E}_{g,n}(\Box, \dots, \Box))$ being a linear combination of products of divisors. 

The rules for multiplying cohomology classes represented by marked stable graphs are explained in~\cite{GraPan}. 
Following these rules, we obtain the product in what is called the {\em strata algebra}. It is well-known that the Picard group of $\oM_{g,n}$ is spanned by $\kappa_1$, $\psi_i$ and the boundary divisors. For $g \geq 3$ these classes form its basis. Denote by $\delta_0$ the boundary divisor of stable curves with a nonseparating node. Consider a product of two divisors and denote by $a$ and $b$ the coefficients of $\delta_0$ in the basis described above. Then the coefficient of every 2-loop in their product is equal to $ab$. Indeed, contracting either edge of a 2-loop we always get $\delta_0$, thus the only way to get a 2-loop as a product of two elements of the basis is when both elements are equal to~$\delta_0$. In particular, the coefficients of all 2-loops in any linear combination of products of divisors are equal to each other, which is not the case in our expression for $\ch_2(\mathbb{E}_{g,n}(\Box, \dots, \Box))$.

Before concluding a contradiction, we must consider whether two different expressions in the strata algebra can represent the same cohomology class. This is the case if and only if their difference is a tautological relation. Thus we are left with the question: is there a tautological relation in the degree 2 part of the strata algebra such that its coefficients have one value for the even 2-loops and another value for the odd 2-loops? A family of tautological relations that conjecturally spans all existing relations was described in~\cite{Pixton} (see also~\cite{PaPiZv}).
While the completeness of the relations is
open in general, specific cases can be verified by
computer calculations (if appropriate pairings are
nondegenerate). 

We have checked that the known relations
indeed span all tautological relations in $RH^4(\oM_{6,2})$. {The
rank of $RH^4(\oM_{6,2})$ is equal to 154: the known relations
provide an upper bound of 154 for the rank, and the pairing
with $RH^{30}(\oM_{6,2})$ provides a lower bound of 154.}

The known relations are given by explicit expressions of degree at least $(g+1)/3$ and by push-forwards of these expressions under gluing maps. Thus we have to consider three cases. First, the push-forwards under gluing maps 
$$
\oM_{g_1, n_1+1} \times \oM_{g_2, n_2+1} \to \oM_{g,n}
$$
do not involve 2-loops at all. Second, the push-forwards under the gluing map
$$
\oM_{g-1,n+2} \to \oM_{g,n}
$$
have degree at least $(g-1+1)/3+1 > 2$ if $g \geq 4$. Thus they don't appear in degree~2. Finally, the expressions on $\oM_{g,n}$ itself have degree at least $(g+1)/3$. Thus for $g \geq 6$ and $n=2$ there are no relations at all in 
degree~2 which involve 2-loops.

To sum up: the coefficients of 2-loops are equal in any linear combination of products of divisors; they are different in our expression for the second Chern character; and the difference cannot be accounted for by tautological relations, at least in the case of $\oM_{6,2}$. \qed

\section {The Verlinde bundle for $\mathfrak {sl}_r$ at level~1}

We consider here the Verlinde bundle for $\mathfrak {sl}_r$ at level $1$. We show that its Chern character has a particularly nice expression over the moduli space $\cM_{g,n}^c$ of curves with compact Jacobians, that is curves with no non-separating nodes; see Proposition \ref{exponential}. 

There are $r$ level $1$  representations denoted $$\omega_0, \ldots, \omega_{r-1},$$ and indexed by Young diagrams with $i$ boxes on a single column, for $0\leq i\leq r-1.$  Note that $$\omega_0^* = \omega_0, \,\, \omega_i^* = \omega_{r-i}, \text{ for }1 \leq i \leq r-1.$$ Moreover,  we have
$$
\mathsf w_i:=\mathsf w (\omega_i) =\frac{i(r-i)}{2r}.
$$ 

The Verlinde bundle $\mathbb E_g(\mu_1, \ldots, \mu_n)\to \cM_{g, n}^c$ depends on the choice of $n$ representations $\mu_1, \ldots, \mu_n$ from the $r$ level $1$ representations $$\omega_0, \ldots, \omega_{r-1}.$$ We will think of the representation $\mu_i$ as an integer in the set $\{0, 1, \ldots, r-1\}$ given by the number of boxes in the column of the corresponding Young diagram; the notation $\mu_i$ will be used for both the representation and for the associated integer. 

The following well-known result simplifies the formulas. The proof is via direct computation using Corollary $9.8$ of \cite {beauville}.

\begin {lemma}\label{vn} For any representations $\mu_{1}, \ldots, \mu_{n}$ at level $1$, we have 
$$
d_g(\mu_1, \ldots, \mu_n)=\begin{cases}
r^g & \mbox{ if } \quad \sum \mu_i \equiv 0 \mod r, \\
0 & \mbox{ otherwise}.
\end{cases}
$$ 
\end {lemma}

Assume that the total number of boxes $\sum \mu_i$ in $\mu_{1}, \ldots, \mu_{n}$ is divisible by $r$. From Theorem 1, we obtain that over $\mathcal M_{g, n}^c$ 
$$
\ch(\mathbb E_g(\mu_1, \ldots, \mu_n))=r^g\exp\left(-\frac{r-1}{2}\lambda_1\right)\cdot$$ $$\sum_{\Gamma} \frac{1}{|\text{Aut}(\Gamma)|} (\iota_{\Gamma})_{\star} \left(\prod_{e} \frac{1-\exp(-\mathsf w_e (\psi'_e+\psi''_e))}{\psi'_e+\psi''_e}\cdot \prod_{i=1}^n e^{-\mathsf w_{\mu_i} \psi_i}\right).
$$ 
The sum is taken over all stable trees $\Gamma$, with the two products  indexed by the edges and legs of the tree $\Gamma$. The weight $\mathsf w_e$ in the edge contribution is determined in a unique way by assigning a remainder modulo~$r$ to each half-edge in such a way that  
\begin {itemize}
\item [(i)] the remainders on the two halves of an edge add up to 0 modulo~$r$;
\item [(ii)] for each vertex~$v$, the remainders on the half-edges and legs adjacent to~$v$ add up to 0 modulo~$r$.
\end {itemize}
Uniqueness can be seen by assigning remainders inductively, starting with a terminal vertex of the tree. 
Since $\mathsf w_i=\mathsf w_{r-i}$, the weight of an edge can be calculated from the remainder of either of its half-edges, by condition~(i) above. 

We show that the expression for the Chern character can be further simplified: 
\begin {proposition}\label{exponential} The Chern character of the $\mathfrak {sl}_r$ level $1$ Verlinde bundle over $\cM_{g, n}^c$ equals 
$$\ch (\mathbb E_g(\mu_1, \ldots, \mu_n))=r^g \exp \left(-\frac{r-1}{2} \lambda_1-\sum_{i=1}^n \mathsf w_{\mu_i} \psi_i +\sum_{e} \mathsf w_e \delta_e\right).$$\end{proposition} The last sum is indexed by the single-edge graphs $e$ (the self-edge excluded), with $\delta_e$ denoting the corresponding boundary divisor.

\paragraph{Proof.} The expression in the Proposition is a consequence of the expansion over $\mathcal M_{g, n}^c$: $$\exp \left(\sum_{e} \mathsf w_e \delta_e\right)=\sum_{\Gamma} \frac{1}{|\text{Aut}(\Gamma)|} (\iota_{\Gamma})_{\star} \left(\prod_{e} \frac{1-\exp(-\mathsf w_e(\psi'_e+\psi''_e))}{\psi_e'+\psi_e''}\right).$$ The equation above can be checked in each fixed codimension $k$ by expanding the two sides and matching terms. Before doing so, let us remark that the graphs appearing on both sides must have all terminal vertices decorated by legs. For otherwise, such terminal vertices would be assigned the integer $0$ by the balancing condition (ii), and then their half-edge contribution $\mathsf w_0$ would vanish. As a consequence, such graphs do not have nontrivial automorphisms. Now, the left hand side can be written as a sum of terms $$\frac{1}{(k_1+1)!\ldots (k_{\ell}+1)!} \mathsf w_{e_1}^{k_1+1} \cdots \mathsf w_{e_{\ell}}^{k_{\ell}+1}\delta_{e_1}^{k_1+1}\cdots \delta_{e_\ell}^{k_\ell+1},$$ where $k_1+\ldots+k_\ell=k-\ell$, and $e_1, \ldots, e_{\ell}$ are one-edge graphs. On the other hand, a graph $\Gamma$ with $\ell$ edges $e_1, \ldots, e_{\ell}$ appearing on the right hand side yields a codimension $k$ contribution given by $$\sum_{k_1+ \ldots+ k_{\ell}=k-\ell} \frac{\mathsf w_{e_1}^{k_1+1}}{(k_1+1)!} \cdots \frac{\mathsf w_{e_\ell}^{k_{\ell}+1}}{(k_{\ell}+1)!} (\iota_{\Gamma})_{\star}\left((-\psi'_{e_1}-\psi''_{e_1})^{k_1}\cdots (-\psi'_{e_\ell}-\psi''_{e_{\ell}})^{k_{\ell}}\right).$$ To match terms, note first that for each choice of edges $e_1, \ldots, e_{\ell}$ which appear in the left hand side expression (giving a non-zero term), there is {\it exactly one} graph $\Gamma$ on the right which yields precisely these edges when all the other edges of the graph are collapsed. The proof is completed noting that $$\delta_{e_j}^{k_j+1}=\iota_{\star} (-\psi'_{e_j}-\psi''_{e_j})^{k_j}.$$ As explained in Proposition \ref{Prop:CompType}, this may fail for boundary divisors of compact type which carry one vertex without markings; however, such divisors were remarked not to appear in the formulas above. \qed

\begin {remark} As a further simplification, in the absence of markings, all half-edges are assigned the integer $0$, by uniqueness. The above formula becomes $$\text{ch} (\mathbb E_g)=r^g\exp\left(-\frac{r-1}{2} \lambda_1\right)\in H^{\star}(\mathcal M_g^c).$$ Up to taking duals, the same answer is found in \cite {vdG} for the Chern character of the projectively flat bundle of level $r$ theta functions over the fibers of the universal Jacobian $$\mathcal J\to \mathcal M_{g}^c.$$ This agreement should be a consequence of an extended strange duality isomorphism,
yet unproven, between level $1$ rank $r$ conformal blocks for compact type curves and level $r$ classical theta functions.  
\end{remark}

\end{document}